\documentclass[12pt,leqno]{article}
\usepackage{latexsym}
\usepackage{bbm}
\usepackage{amssymb,amsmath}
\begin{titlepage}
\title{On the large sieve with sparse sets of moduli}
\author{Stephan Baier}
\date{22.08.05}
\end{titlepage}
\begin{document}
\maketitle
$ $\\
$ $\\
{\bf Address of the author:}\medskip\\ 
Stephan Baier\\  
Jeffery Hall\\ 
Department of Mathematics and Statistics\\
Queen's University\\
University Ave\\
Kingston, Ontario, Canada\\
K7L 3N6\medskip\\
e-mail: sbaier@mast.queensu.ca
\newpage 
$ $\\
{\bf Abstract:} Extending a method of D. Wolke \cite{Wol}, 
we establish a general 
result on the large sieve with sparse
sets ${\cal{S}}$ of moduli which are in a sense well-distributed in 
arithmetic 
progressions. We then use this result together with Fourier techniques to
obtain large sieve bounds for the case when ${\cal{S}}$ consists
of sqares. These bounds improve a recent result by L. Zhao \cite{Zha}.
\section{A general result on the large sieve}
Throughout this paper, we reserve the symbols 
$c_i$ $(i=1,2,...)$ for absolute 
constants and the symbol $\varepsilon$ for an arbitrary (small) positive number. The $\ll$-constants in our estimates may depend on $\varepsilon$. As usual in analytic number theory, the $\varepsilon$ may be different from line to line. We further suppose that 
$(a_n)$ is a sequence of complex numbers and that $Q,N\ge 1$. We set
\begin{equation}
S(\alpha):=\sum\limits_{n\le N} a_n e(n\alpha).\label{0}
\end{equation}

Bombieri's \cite{Cla} classical large sieve inequality asserts that 
\begin{equation}
\sum\limits_{q\le Q} \sum\limits_{\scriptsize \begin{array}{cccc} 
a=1\\(a,q)=1\end{array}}^q \left\vert S\left(\frac{a}{q}\right)\right\vert^2
\le (N+Q^2)Z,\label{1}
\end{equation}
where 
$$
Z:=\sum\limits_{n\le N} \vert a_n\vert^2.
$$
One may ask whether (\ref{1}) can be improved if 
the moduli $q$ run over a sparse set ${\cal{S}}$ of natural numbers $\le Q$.
It seems to be
difficult to obtain a considerable improvement if nothing is known about
the structure of ${\cal{S}}$. The goal of the present paper is to improve 
(\ref{1}) for sets ${\cal{S}}$ of moduli which are in a sense well-distributed 
in arithmetic progressions.  

In the sequel, we suppose, more generally, that ${\cal{S}}\subset (M,M+Q]$,
where $0\le M\le Q$. We put $S:=\left\vert {\cal{S}}\right\vert$ (the cardinality of ${\cal{S}}$). 
For $t\in \mathbbm{N}$ we put
$$
{\cal{S}}_t:=\{q\in \mathbbm{N}\ :\ tq\in {\cal{S}}\}
$$
and $S_t:=\vert {\cal{S}}_t \vert$. 
We note that ${\cal{S}}_t\subset 
(M/t,M/t+Q/t]$. We shall require that the number of 
elements of ${\cal{S}}_t$ in short segments of arithmetic progressions
does not differ too much from the expected number. To measure the 
distribution of ${\cal{S}}_t$ in segments of arithmetic progressions,
we define the quantity 
$$
A_t(u,k,l):=\max\limits_{M/t\le y\le (M+Q)/t} 
\vert \{q\in {\cal{S}}_t\cap (y,y+u] \ : \ q\equiv l\mbox{ mod }k\} \vert,
$$
where $u\ge 0$, $k\in \mathbbm{N}$, $l\in \mathbbm{Z}$ with $(k,l)=1$. 
We shall establish the following\\

{\bf Theorem 1:} \begin{it} We have
\begin{eqnarray*}
\label{4} & &\\
\sum\limits_{q\in {\cal{S}}} \sum\limits_{\scriptsize \begin{array}{cccc} 
a=1\\(a,q)=1\end{array}}^q \left\vert S\left(\frac{a}{q}\right)\right\vert^2
&\le& c_1NZ\left(1+\max\limits_{r\le \sqrt{N}}\ 
\max\limits_{1/N\le z\le 1/(r\sqrt{N})} \max\limits_{\scriptsize
\begin{array}{cccc} h\in \mathbbm{Z}\\(h,r)=1\end{array}}\right. \nonumber\\
& & \left. \sum\limits_{t\vert r}
\sum\limits_{\scriptsize \begin{array}{cccc} 
0< \vert m\vert\le 6rzQ/t\\ (m,r/t)=1\end{array}} 
A_t\left(\frac{2Q}{tzN},\frac{r}{t},hm\right)\right).\nonumber
\end{eqnarray*}
\end{it}

If we assume 
the set ${\cal{S}}_t$ to be nearly evenly distributed in the residue classes 
$l$ mod $k$, then, if $M/t\le y\le (M+Q)/t-u$, the expected number of elements of the set 
$$
\{q\in {\cal{S}}_t\cap (y,y+u] \ : \ q\equiv l\mbox{ mod }k\} 
$$
is
$$
\frac{S_t/k}{Q/t}\cdot u.
$$
Therefore, if ${\cal{S}}_t$ is well-distributed in the residue classes $l$ mod $k$, we may expect, for any $u\ge 0$, that 
\begin{equation}
A_t(u,k,l)\le \left(1+\frac{S_t/k}{Q/t}\cdot u\right)X, \label{20}
\end{equation}
where $X\ge 1$ is small compared to $Q$ and $N$.

By a short calculation, we infer the following bound from Theorem 1.\\ 

{\bf Theorem 2:} \begin{it}   
Suppose the condition (\ref{20}) to hold for all $t,k,l,u$ with  
$t\le \sqrt{N}$, $k\le \sqrt{N}/t$, $(k,l)=1$ and  
$kQ/\sqrt{N}\le u\le Q/t$. Then
\begin{equation}
\sum\limits_{q\in {\cal{S}}} \sum\limits_{\scriptsize \begin{array}{cccc} 
a=1\\(a,q)=1\end{array}}^q \left\vert S\left(\frac{a}{q}\right)\right\vert^2
\le c_2\left(N+ QXN^{\varepsilon}\left(\sqrt{N}+S\right)\right)Z. 
\label{8}
\end{equation} 
\end{it} 

This is stronger than the classical large sieve inequality
(\ref{1}) if $N^{1+\varepsilon}\ll 
QX(\sqrt{N}+S)\ll Q^{2-\varepsilon}$. If the set ${\cal{S}}$ is really sparse, that is,
if $S$ is small compared to $Q$, and if the condition (\ref{20}) holds with $X=N^{\varepsilon}$, then (\ref{8}) is sharper than (\ref{1}) if $Q\gg N^{1/2+\varepsilon}$. In section 6 we shall see that in the case of square moduli $X=N^{\varepsilon}$ is an admissible choice
in (\ref{20}).

A conjecture of Elliott [5] would imply that the left-hand side of (\ref{8}) is bounded by
\begin{equation}
\ll (N+QS)Z \label{El}
\end{equation}
if ${\cal{S}}$ contains only primes. From (\ref{8}), we obtain
the slightly weaker bound $\ll (N+QN^{\varepsilon}S)Z$ if
$X=N^{\varepsilon}$ is admissible and $S\gg \sqrt{N}$. 

\section{The case of squares}
Recently, L. Zhao \cite{Zha} studied the case when the moduli
$q$ are squares, that is, he investigated the order of magnitude of the expression 
$$
T:=
\sum\limits_{q\le Q} \sum\limits_{\scriptsize \begin{array}{cccc} a=1\\(a,q)=1\end{array}}^{q^2}
\left\vert S\left(\frac{a}{q^2}\right)\right\vert^2.
$$
He proved the estimate
\begin{equation}
T\ll (\log 2Q)\left(Q^3+(N\sqrt{Q}+\sqrt{N}Q^2)N^{\varepsilon}\right)Z
\label{3}
\end{equation}
and conjectured that 
\begin{equation}
T\ll Q^{\varepsilon}(Q^3+N)Z.\label{LZ}
\end{equation}
The classical form (\ref{1}) of the large sieve implies only the bound
\begin{equation}
T\ll (N+Q^{4})Z, \label{Z1}
\end{equation}
which is weaker than (\ref{3}) if $Q\gg N^{2/7+\varepsilon}$. 
Using the bound
\begin{equation}
\sum\limits_{\scriptsize \begin{array}{cccc} 
a=1\\(a,q)=1\end{array}}^{q^2}
\left\vert S\left(\frac{a}{q^2}\right)\right\vert^2
\ll (N+q^2)Z, \label{ind}
\end{equation}
which follows from our later
Lemma 1 with $\Delta=1/q^2$ and $(\alpha_r)$ beeing the sequence formed by all fractions $a/q^2$ with $1\le a\le q^2$
and $(a,q)=1$, we also obtain the bound 
\begin{equation}
T \ll Q(N+Q^{2})Z \label{Z2}
\end{equation}
by summing up (\ref{ind}) over all $q\le Q$. This bound 
is weaker than (\ref{3}) if $Q\ll N^{1/2-\varepsilon}$. 
Thus, (\ref{3}) is sharper
than both (\ref{Z1}) and (\ref{Z2}) if $N^{2/7+\varepsilon}\ll Q\ll 
N^{1/2-\varepsilon}$. 
   
Employing Theorem 2 with ${\cal{S}}$ a set of squares, we shall 
obtain the following improvement of Zhao's bound
(\ref{3}).\\

{\bf Theorem 3:} \begin{it} We have
\begin{equation}
T \ll (\log 2Q)N^{\varepsilon}(Q^{3}+N+N^{1/2}Q^2)Z. 
\label{9}
\end{equation}\end{it}

The bound (\ref{9}) is sharper than the three bounds (\ref{3}), (\ref{Z1})
and (\ref{Z2}) if 
$N^{1/4+\varepsilon}\ll Q\ll N^{1/3-\varepsilon}$. Combining the elementary
methods which we will use for the proof of Theorem 3 with Fourier analytic techniques,
we shall further prove\\

{\bf Theorem 4:} \begin{it} We have
$$
T\ll \left\{ \begin{array}{llll} Q^{3/5+\varepsilon}NZ, &
\mbox{ if } Q\le N^{5/12},\\ \\ Q^{3+\varepsilon}Z, &
\mbox{ if } Q> N^{5/12}.\end{array} \right.
$$
\end{it} 

This bound is sharper than (\ref{3}), (\ref{Z1}) and (\ref{Z2}) if 
$N^{5/14+\varepsilon}\ll Q \ll N^{1/2-\varepsilon}$.
Moreover, it establishes Zhao's conjecture (\ref{LZ}) for $Q\gg N^{5/12}$.

\section{The case of primes}     
For the case when ${\cal{S}}$ is the full set of all primes 
$p\le Q$ D. Wolke \cite{Wol} proved the estimate 
\begin{equation}
\sum\limits_{p\le Q} \sum\limits_{a=1}^{p-1} 
\left\vert S\left(\frac{a}{p}\right)\right\vert^2
\le \frac{c_3}{1-\delta} 
\frac{Q^2\log\log Q}{\log Q}Z \label{2}
\end{equation}
provided that 
\begin{equation}
Q\ge 10,\ \ \ \ \ \  N=Q^{1+\delta},\ \ \ \ \ \  0<\delta<1.\label{QCON}
\end{equation} 
In this range Elliott's conjecture (\ref{El}) would give the slightly better bound $\ll Q^2Z/\log Q$. 

Now we want to prove that Theorem 1 with $M=0$ and ${\cal{S}}$ beeing the set of all primes $p\le Q$ implies Wolke's bound (\ref{2}). 
We need to estimate the term $A_t(u,k,l)$. First we consider the case when $t=1$. By the Brun-Titchmarsh inequality, we have
\begin{equation}
A_1\left(\frac{2Q}{zN},r,l\right)\le 
\frac{4Q}{zN\varphi(r)\log (2Q/rzN)}\label{BT1}
\end{equation}
if $2Q/(zN)>r$. If $rz\le 1/\sqrt{N}$, then $2Q/(zN)>r$ is satisfied since $1/\sqrt{N}< 2Q/N$ by (\ref{QCON}). 
From (\ref{QCON}) and (\ref{BT1}), we deduce
\begin{equation}
\sum\limits_{\scriptsize \begin{array}{cccc} 
0< \vert m\vert\le 6rzQ\\ (m,r)=1\end{array}} 
A_1\left(\frac{2Q}{zN},r,hm\right)\le 
c_4\frac{Q^2\log\log Q}{N(1-\delta)\log Q} \label{BT2}
\end{equation}
for any integer $h$ with $(r,h)=1$.

If $t\ge 2$, then ${\cal{S}}_t$ contains at most
1 element. This implies 
\begin{eqnarray}
\label{BT3} & &
\sum\limits_{\scriptsize \begin{array}{cccc} 
t\vert r\\ t\ge 2\end{array}} 
\sum\limits_{\scriptsize \begin{array}{cccc} 
0< \vert m\vert\le 6rzQ/t\\ (m,r/t)=1\end{array}} 
A_t\left(\frac{2Q}{tzN},\frac{r}{t},hm\right)\\
&\le& 
\sum\limits_{t\vert r} \frac{12rzQ}{t} \le 
c_{5}rzQ\log\log 10r \le 
c_{6}\frac{Q\log\log Q}{\sqrt{N}}\nonumber
\end{eqnarray}
if $r\le \sqrt{N}< Q$ and $z\le 1/(r\sqrt{N})$. 
Using $N=Q^{1+\delta}$,
it is easy to check that there exists a constant $c_{7}$ such that we have
\begin{equation}
\frac{Q\log\log Q}{\sqrt{N}}\le 
c_{7}\frac{Q^2\log\log Q}{N(1-\delta)\log Q}\label{BT4}
\end{equation} 
for all $Q\ge 10$ and $0<\delta<1$. From Theorem 1, (\ref{BT2}),
(\ref{BT3}) and (\ref{BT4}), we obtain Wolke's bound (\ref{2}).  

\section{Counting Farey fractions in short intervals}
In this section we establish some preliminary results which we then use for the proof of Theorem 1. 
Our starting point is the following variant of the large
sieve which follows immediately from Theorem 2.11 in \cite{Lem}.\\ 

{\bf Lemma 1:} \begin{it} Let $\left(\alpha_r\right)_{r\in\mathbbm{N}}$ be a
sequence of real numbers. Suppose that $0<\Delta\le 1/2$ and 
$R\in \mathbbm{N}$. Put 
$$
K(\Delta):=\max\limits_{\alpha\in \mathbbm{R}} 
\sum\limits_{\scriptsize \begin{array}{cccc} r=1\\
\vert\vert \alpha_r -\alpha\vert\vert\le \Delta \end{array}}^R 1,
$$
where $\vert\vert x \vert\vert$ denotes the distance of a real $x$
to its closest integer.
Then 
$$
\sum\limits_{r=1}^R \left\vert S\left(\alpha_r\right)\right\vert^2
\le c_{8}K(\Delta)(N+\Delta^{-1})Z.
$$
\end{it}\\

In our situation, the sequence $\alpha_1,...,\alpha_R$ equals the sequence
of Farey fractions $a/q$ with $q\in {\cal{S}}$, $1\le a\le q$ and $(a,q)=1$. 
For $\alpha\in \mathbbm{R}$ we put
$$
I(\alpha):=[\alpha-\Delta,\alpha+\Delta]
\ \ \mbox{ and }\ \ 
P(\alpha):= \sum\limits_{\scriptsize \begin{array}{cccc} 
q\in {\cal{S}}, (a,q)=1\\
a/q\in I(\alpha) \end{array}} 1.
$$
Then we have
\begin{equation}
K(\Delta)=\max\limits_{\alpha\in \mathbbm{R}} P(\alpha).\label{60}
\end{equation}

To estimate $P(\alpha)$, we begin with a method of D. Wolke \cite{Wol}. Let
\begin{equation}
\tau:=\frac{1}{\sqrt{\Delta}}.\label{P1}
\end{equation}
Then, by Dirichlet's approximation theorem, $\alpha$ can be written in the form
\begin{equation}
\alpha=\frac{b}{r}+z, \ \ \mbox{ where }\ \  r\le \tau,\ (b,r)=1,\ 
\vert z\vert \le \frac{1}{r\tau}.\label{P2}
\end{equation}
Thus, it suffices to estimate $P(b/r+z)$ for all $b,r,z$ satisfying
(\ref{P2}).

We further note that we can restrict ourselves to the case when 
\begin{equation}
z\ge \Delta.\label{P4}
\end{equation}
If $\vert z\vert<\Delta$, then
$$ 
P(\alpha)\le P\left(\frac{b}{r}-\Delta\right)+P\left(\frac{b}{r}+\Delta\right).
$$ 
Furthermore, we have
$$
\Delta=\frac{1}{\tau^2}\le \frac{1}{r\tau}.
$$ 
Therefore
this case can
be reduced to the case $\vert z\vert=\Delta$. 
Moreover, as $P(\alpha)=P(-\alpha)$, we can choose
$z$ positive. So we can assume (\ref{P4}).   

Summarizing the above observations, we deduce\\

{\bf Lemma 2:} \begin{it} We have  
\begin{equation}
K(\Delta)\le 2\max\limits_{\scriptsize \begin{array}{cccc} 
r\in \mathbbm{N} \\ r\le 1/\sqrt{\Delta} \end{array}} 
\max\limits_{\scriptsize \begin{array}{cccc} b\in \mathbbm{Z}\\ 
(b,r)=1\end{array}}
\max\limits_{\Delta\le z\le \sqrt{\Delta}/r} P\left(\frac{b}{r}+z\right).
\end{equation} \end{it}

The next lemma provides a first estimate for $P\left(b/r+z\right)$.\\
       
{\bf Lemma 3:} \begin{it}
Suppose that the conditions (\ref{P1}), (\ref{P2}) and (\ref{P4}) are 
satisfied. Suppose further that 
\begin{equation}
\frac{Q\Delta}{z}\le \delta\le Q. \label{P10}
\end{equation}
Let $I(\delta,y):=[y-\delta,y+\delta]$, $J(\delta,y):=
[(y-4\delta)rz,(y+4\delta)rz]$ and 
$$
\Pi(\delta,y):=\sum\limits_{q \in {\cal{S}}\cap I(\delta,y)}
\sum\limits_{\scriptsize 
\begin{array}{cccc} m\in J(\delta,y) \smallskip
\\ m \equiv -bq \mbox{ mod } r\smallskip \\ m\not= 0
\end{array}} 1.
$$
Then,
$$
P\left(\frac{b}{r}+z\right) \le 2+
\frac{1}{\delta} \int\limits_{M}^{M+Q} \Pi(\delta,y)\ \rm{d}y.
$$
\end{it}\\

{\bf Proof:} By $\delta\le Q$, we have 
\begin{equation}
P(\alpha)\le \frac{1}{\delta} \int\limits_{M}^{M+Q} P(\alpha,y,\delta) 
\ {\rm d} y,\label{P6}
\end{equation}
where  
$$
P(\alpha,y,\delta):=\sum\limits_{\scriptsize \begin{array}{cccc} 
q \in {\cal{S}}\cap I(\delta,y) \smallskip\\ (a,q)=1\smallskip\\
a/q\in I(\alpha) \end{array}} 1.
$$
Now, for  $a/q\in I(\alpha)$,
we have 
$$
q(\alpha-\Delta)\le a \le q(\alpha+\Delta).
$$ 
From this and $\alpha=b/r+z$, we obtain
\begin{equation}
qr(z-\Delta)\le ar-bq\le qr(z+\Delta). \label{P20}
\end{equation}
If $y-\delta\le q\le y+\delta$, then from (\ref{P2}), (\ref{P4}), (\ref{P10}) and (\ref{P20}) it follows that
\begin{equation}
(y-4\delta)rz\le (y-\delta)r(z-\Delta)\le ar-bq\le 
(y+\delta)r(z+\Delta)\le (y+4\delta)rz.\label{P21}
\end{equation}  
If $ar-bq=0$, then $r=q$ because $(a,q)=1=(b,r)$. From this
observation, (\ref{P6}) and (\ref{P21}), we deduce the result
of Lemma 3.$\Box$

\section{Proof of Theorem 1}

Next, we express $\Pi(y,\delta)$ in terms of $A_t(u,k,l)$. This shall lead
us to the following estimate for $P(b/r+z)$.\\

{\bf Lemma 4:} \begin{it} We have
$$
P\left(\frac{b}{r}+z\right) \le 2+c_9
\sum\limits_{t\vert r} \sum\limits_{\scriptsize \begin{array}{cccc} 
0<\vert m\vert \le 6rzQ/t\\ (m,r/t)=1\end{array}} 
A_t\left(\frac{2\Delta Q}{tz},\frac{r}{t},-\overline{b}m\right),
$$
where $\overline{b}b$ $\equiv$ $1$ mod $r$.
\end{it}\\

On choosing $\Delta:=1/N$, Theorem 1 follows immediately from
Lemmas 1,2 and 4.\\ 

{\bf Proof of Lemma 4:} 
We split $\Pi(\delta,y)$ into
$$
\Pi(\delta,y)=\sum\limits_{t\vert r} 
\sum\limits_{\scriptsize 
\begin{array}{cccc} q \in {{\cal{S}}_t}\cap I(\delta/t,y/t)\smallskip\\
(q,r/t)=1 \end{array}}
\sum\limits_{\scriptsize 
\begin{array}{cccc} m\in J(\delta/t,y/t) \smallskip
\\ m \equiv -bq \mbox{ mod } r/t\smallskip \\ m\not= 0
\end{array}} 1. 
$$
Rearranging the order of summation, and using the definition of 
$A_t(u,k,l)$, the right-hand side is
\begin{eqnarray}
\label{AAA} &=&\sum\limits_{t\vert r} 
\sum\limits_{\scriptsize 
\begin{array}{cccc} m\in J(\delta/t,y/t) \smallskip
\\ (m,r/t)=1\smallskip \\ m\not= 0
\end{array}}
\sum\limits_{\scriptsize 
\begin{array}{cccc} q \in {{\cal{S}}_t}\cap I(\delta/t,y/t)\smallskip\\
q\equiv -\overline{b}m \mbox{ mod } r/t \end{array}} 1\\
&\le &\sum\limits_{t\vert r} 
\sum\limits_{\scriptsize 
\begin{array}{cccc} m\in J(\delta/t,y/t) \smallskip
\\ (m,r/t)=1\smallskip \\ m\not= 0 \end{array}}
A_t\left(\frac{2\delta}{t},\frac{r}{t},-\overline{b}m\right).
\nonumber
\end{eqnarray}
Integrating the last line of (\ref{AAA}) over $y$ in the interval $M\le y\le M+Q$, and rearranging
the order of summation and integration, we obtain

\begin{equation}
\int\limits_{M}^{M+Q} \Pi(\delta,y)\ {\rm d}y
\le 2\delta \sum\limits_{t\vert r} 
\sum\limits_{\scriptsize 
\begin{array}{cccc} (M-4\delta)rz/t\le m\le (M+Q+4\delta)rz/t \smallskip
\\ (m,r/t)=1\smallskip \\ m\not= 0 \end{array}}
A_t\left(\frac{2\delta}{t},\frac{r}{t},-\overline{b}m\right). \label{L4}
\end{equation}
Choosing $\delta:=Q\Delta/z$, and taking $0\le M\le Q$ and $Q\Delta/z\le Q$
into account, we obtain the result of Lemma 4 from 
Lemma 3 and (\ref{L4}).
$\Box$\\

From Lemma 4, we also infer the following estimate for 
$P(b/r+z)$ by a short calculation.\\    

{\bf Lemma 5:} \begin{it} Suppose that the conditions (\ref{P1}), (\ref{P2}) 
and (\ref{P4}) are satisfied.  
Suppose further the condition (\ref{20}) to hold for   
$t\vert r$, $k=r/t$, $(k,l)=1$ and  
$u=2\Delta Q/(tz)$. Then
$$
P\left(\frac{b}{r}+z\right)\le c_{10}\left(1+ QX
\Delta^{-\varepsilon}\left(rz+\Delta S\right)\right)Z. 
$$\end{it}

This estimate corresponds to Theorem 2. We shall use it in section 7.

\section{Proof of Theorem 3}
In this section, we derive Theorem 3 from Theorem 2.
First, we rewrite the sum in question in the form
$$
T=\sum\limits_{q\le Q} \sum\limits_{\scriptsize \begin{array}{cccc} 
a=1\\(a,q)=1\end{array}}^{q^2}
\left\vert S\left(\frac{a}{q^2}\right)\right\vert^2
= \sum\limits_{q\in {\cal{S}}} \sum\limits_{\scriptsize \begin{array}{cccc} 
a=1\\(a,q)=1\end{array}}^q
\left\vert S\left(\frac{a}{q}\right)\right\vert^2, 
$$
where ${\cal{S}}$ is the set of squares $\le Q^2$. 
We split up the set ${\cal{S}}$ into 
$O(\log Q)$ subsets of the form
$$
{\cal{S}}(Q_0):={\cal{S}}\cap (Q_0,2Q_0],
$$
where $Q_0\ge 1$.
Our aim is to estimate the corresponding partial sums.
As previously, we define 
$$
{\cal{S}}_t(Q_0):=\{q\in \mathbbm{N}\ :\ tq\in 
{\cal{S}}(Q_0)\}
$$
and $S_t(Q_0):=\vert {\cal{S}}_t(Q_0) \vert$. 
We now determine the set ${\cal{S}}_t(Q_0)$. Let $t=p_1^{v_1}\cdots p_n^{v_n}$
be the prime number factorization of $t$. For $i=1,...,n$ let
$$
u_i:=\left\{\begin{array}{llll} v_i, & \mbox{ if } v_i \mbox{ is even,}\\ \\
v_i+1, & \mbox{ if } v_i \mbox{ is odd.} \end{array}\right.
$$
Put 
$$
f_t:=p_1^{u_1/2}\cdots p_n^{u_n/2}.
$$
Then $q=q_1^2\in{\cal{S}}$ is divisible by $t$ iff $q_1$ is divisible by
$f_t$. Thus, 
$$
{\cal{S}}_t(Q_0)=\left\{q_2^2g_t\ : \ \sqrt{Q_0}/f_t<q_2\le \sqrt{2Q_0}/f_t\right\}
\subset (Q_0/t,2Q_0/t],
$$
where  
$$
g_t:=\frac{f_t^2}{t}=p_1^{u_1-v_1}\cdots p_n^{u_n-v_n}.
$$
As previously, we suppose that $u\ge 0$,
$k\in \mathbbm{N}$, $l\in \mathbbm{Z}$ and $(k,l)=1$, and 
define 
$$
A_t(u,k,l):=\max\limits_{Q_0/t\le y\le 2Q_0/t} 
\vert \{q\in {\cal{S}}_t(Q_0)\cap (y,y+u] \ : \ q\equiv l\mbox{ mod }k\} \vert.
$$
Let $\delta_t(k,l)$ be the number of solutions $x$ mod $k$ to the congruence
\begin{equation}
x^2g_t\equiv l \mbox{ mod } k. \label{23}
\end{equation}
Then, from our above observations it follows quickly that 
$$
A_t(u,k,l)\le c_{11}
\left(1+\frac{S_t/k}{Q/t}\cdot u\right)\delta_t(k,l). 
$$
The remaining task is to bound $\delta_t(k,l)$.

If $(g_t,k)>1$, then $\delta_t(k,l)=0$ since $k$ and $l$ are supposed to be
coprime. Therefore, we can assume that $(g_t,k)=1$. Let $g$ mod $k$ 
be the multiplicative inverse of
$g_t$ mod $k$, i.e. $gg_t\equiv 1$ mod $k$. Put $l^*=gl$. 
Then $(\ref{23})$ is 
equivalent to
\begin{equation}
x^2\equiv l^* \mbox{ mod } k. \label{24}
\end{equation} 
Taking into account that $(k,l^*)=1$, and using some elementary facts on
the number 
of solutions of polynomial congruences modulo prime powers (see \cite{Qua}, 
for example), we see that  
(\ref{24}) has at most $2$ solutions if $k$ is a power of an odd prime and
at most $4$ solutions if $k$ is a power of 2. From this it follows that
for all $k\in \mathbbm{N}$ we have
$$
\delta_t(k,l)\le 2^{\omega(k)+1},
$$   
where $\omega(k)$ is the number of distinct 
prime divisors of $k$. For $k\le \sqrt{N}$
we have 
$$
2^{\omega(k)}\ll N^{\varepsilon}
$$
(see \cite{Han}). Therefore, (\ref{20}) holds with 
\begin{equation}
X:=c_{12}N^{\varepsilon}.\label{26}
\end{equation} 
Now, from Theorem 2, (\ref{26}) and $S\ll \sqrt{Q_0}$, we obtain 
$$
\sum\limits_{q\in {\cal{S}}(Q_0)} 
\sum\limits_{\scriptsize \begin{array}{cccc} 
a=1\\(a,q)=1\end{array}}^q \left\vert S\left(\frac{a}{q}\right)\right\vert^2\le c_{13}\left(N+Q_0N^{\varepsilon}
\left(\sqrt{N}+ \sqrt{Q_0}\right)\right)Z.
$$
This implies the result of Theorem 3. $\Box$

\section{Proof of Theorem 4}
Throughout this section, we suppose that ${\cal{S}}$ consists
of all squares in the interval $(Q_0,2Q_0]$. To prove Theorem
4, we use the following estimates for $P(b/r+z)$.\\

{\bf Lemma 6:} \begin{it} Suppose that the conditions 
(\ref{P1}), (\ref{P2}) and (\ref{P4}) are satisfied. 
Then we have
\begin{equation}
P\left(\frac{b}{r}+z\right)\le c_{14}\Delta^{-\varepsilon}
\left(1+Q_0rz+Q_0^{3/2}\Delta\right)\label{R2}
\end{equation}
and
\begin{equation}
P\left(\frac{b}{r}+z\right)\le 
c_{15}\Delta^{-\varepsilon}\left(Q_0^{3/2}\Delta
+Q_0^{1/2}\Delta r^{-1/2}z^{-1}+\Delta^{-1/4}\right).\label{QQ}
\end{equation}\end{it}

The inequality (\ref{R2}) follows immediately from Lemma 5
and the fact that (\ref{20}) holds with 
$X:=c_{16}\Delta^{-\varepsilon}$ under the conditions of
Lemma 5. 
This can be seen in the same way as it was proved
that (\ref{26}) is an admissable choice in (\ref{20})
under the conditions of Theorem 3. 
 
We postpone the proof of $(\ref{QQ})$ to the last section.

We are now ready to prove Theorem 4. Combining (\ref{R2}) and (\ref{QQ}), we obtain

\begin{equation}
P\left(\frac{b}{r}+z\right)\le 
c_{17}\Delta^{-\varepsilon}\left(Q_0^{3/2}\Delta+
\min\left\{Q_0rz,Q_0^{1/2}\Delta r^{-1/2}z^{-1}\right\}+
\Delta^{-1/4}\right).\label{E1}
\end{equation}
If
$$
z\le \Delta^{1/2}Q_0^{-1/4}r^{-3/4},
$$
then 
$$
\min\left\{Q_0rz,Q_0^{1/2}\Delta r^{-1/2}z^{-1}\right\}=
Q_0rz\le Q_0^{3/4}\Delta^{1/2}r^{1/4}.
$$
If 
$$
z> \Delta^{1/2}Q_0^{-1/4}r^{-3/4},
$$
then 
$$
\min\left\{Q_0rz,Q_0^{1/2}\Delta r^{-1/2}z^{-1}\right\}=
Q_0^{1/2}\Delta r^{-1/2}z^{-1}\le Q_0^{3/4}\Delta^{1/2}r^{1/4}.
$$
From the above inequalities and (\ref{P2}), we deduce that
\begin{equation}
\min\left\{Q_0rz,Q_0^{1/2}\Delta r^{-1/2}z^{-1}\right\}\le 
Q_0^{3/4}\Delta^{3/8}. \label{E2}
\end{equation}
Furthermore,
\begin{equation}
Q_0^{3/4}\Delta^{3/8}= \sqrt{(Q_0^{3/2}\Delta)\cdot \Delta^{-1/4}}
\le Q_0^{3/2}\Delta+\Delta^{-1/4}. \label{E3}
\end{equation}
Combining (\ref{E1}), (\ref{E2}) and (\ref{E3}), we get
\begin{equation}
P\left(\frac{b}{r}+z\right)\le 
c_{18}\Delta^{-\varepsilon}\left(Q_0^{3/2}\Delta+\Delta^{-1/4}\right).
\label{E4}
\end{equation}

Now we choose 
$$
\Delta:=\left\{ \begin{array}{llll} Q_0^{-6/5}, &  \mbox{ if } Q_0\le N^{5/6},\\ \\ N^{-1}, & \mbox{ otherwise.} 
\end{array} \right.
$$ 
Then from Lemma 1, Lemma 2 and (\ref{E4}) it 
follows that 
\begin{equation}
\sum\limits_{\sqrt{Q_0}\le q\le \sqrt{2Q_0}} 
\sum\limits_{\scriptsize \begin{array}{cccc} 
a=1\\(a,q)=1\end{array}}^{q^2}
\left\vert S\left(\frac{a}{q^2}\right)\right\vert^2
\ll \left\{ \begin{array}{llll} Q_0^{3/10+\varepsilon}NZ, &  \mbox{ if } Q_0\le N^{5/6},\\ \\ Q_0^{3/2+\varepsilon}Z, & 
\mbox{ otherwise.} 
\end{array} \right. 
\label{E5}
\end{equation}
We can devide the interval $[1,Q]$ into $O(\log Q)$ subintervals of the
form $\left[\sqrt{Q_0},\sqrt{2Q_0}\right]$, where $1\le Q_0\le Q^2$. Hence,
the result of Theorem 4 follows from (\ref{E5}). $\Box$

\section{Tools from harmonic analysis}
For the proof of (\ref{QQ}) we need the following standard results from harmonic analysis.\\ 

{\bf Lemma 7:} (Poisson summation formula, \cite{Bum}) \begin{it}  
 Let $f(X)$ be a complex-valued 
function on the real numbers that is piecewise continuous with only finitely 
many discontinuities and for all real numbers $a$ satisfies
$$
f(a)=\frac{1}{2}\left(\lim\limits_{x\rightarrow a^-} f(x) +
\lim\limits_{x\rightarrow a^+} f(x)\right).
$$
Moreover, suppose that $f(x)\le c_{19}(1+\vert x\vert)^{-c}$ for some $c>1$.
Then,
$$
\sum\limits_{n\in \mathbbm{Z}} f(n) = \sum\limits_{n\in \mathbbm{Z}} \hat{f}(n),
$$
where 
$$
\hat{f}(x):=\int\limits_{-\infty}^{\infty} f(y)e(xy) {\rm d}y,
$$
the Fourier transform of $f(x)$. \end{it}\\

{\bf Lemma 8:} (see \cite{Zha}, for example) \begin{it}
For $x\in \mathbbm{R}\setminus \{0\}$ define
$$
\phi(x):=\left(\frac{\sin \pi x}{2x}\right)^2.
$$
Set 
$$
\phi(0):=\lim\limits_{x\rightarrow 0}\phi(x)=\frac{\pi^2}{4}.
$$ 
Then $\phi(x)\ge 1$ for $\vert x\vert \le 1/2$, and
the Fourier transform of the function $\phi(x)$ is
$$
\hat{\phi}(s)=\frac{\pi^2}{4}\max\{1-\vert s\vert, 0\}.
$$\end{it}

{\bf Lemma 9:} (see Lemma 3.1. in \cite{Kol}) \begin{it} Let $F$ $:$ 
$[a,b]\rightarrow \mathbbm{R}$ be twice differentiable. 
Assume that $\vert F^{\prime}(x)\vert \ge u>0$ 
for all $x\in [a,b]$. Then,
$$
\left\vert \int\limits_{a}^{b} e^{iF(x)} {\rm d}x\right\vert \ll
\frac{c_{20}}{u}.
$$
\end{it}

{\bf Lemma 10:} (see Lemma 4.3.1. in \cite{Bru}) \begin{it} Let $F$ $:$ 
$[a,b]\rightarrow \mathbbm{R}$ be twice continuously differentiable. 
Assume that $\vert F^{\prime\prime}(x)\vert \ge u>0$ 
for all $x\in [a,b]$. Then,
$$
\left\vert \int\limits_{a}^{b} e^{iF(x)} {\rm d}x\right\vert \le
\frac{c_{21}}{\sqrt{u}}.
$$\end{it}\medskip

We shall also need the following estimate for quadratic Gau\ss{} sums.\\

{\bf Lemma 11:} (see page 93 in \cite{Kol}) \begin{it}
Let $c\in \mathbbm{N}$, $k,l\in 
\mathbbm{Z}$ with $(k,c)=1$. Then,
$$
\sum\limits_{d=1}^r e\left(\frac{kd^2+ld}{c}\right) \le \sqrt{2c}.
$$ 
\end{it}

\section{Proof of (\ref{QQ})} 
Applying Lemma 3 with $Q$ replaced by $Q_0$, $M=Q_0$ and\\
$S:=\left\{q^2 \ : \ \sqrt{Q_0}<q\le \sqrt{2Q_0}\right\}$,
we have, for any $\delta$ satisfying the condition (\ref{P10}), 
\begin{equation}
P\left(\frac{b}{r}+z\right) \le 2+
\frac{1}{\delta} \int\limits_{Q_0}^{2Q_0} 
\Pi(\delta,y)\ \rm{d}y, \label{EIN1}
\end{equation}
where 
$$
\Pi(\delta,y):=\sum\limits_{\sqrt{y-\delta}\le q \le \sqrt{y+\delta}} \sum\limits_{\scriptsize 
\begin{array}{cccc} m\in J(\delta,y) \smallskip
\\ m \equiv -bq^2 \mbox{ mod } r\smallskip \\ m\not= 0
\end{array}} 1.
$$
By Taylors formula and $\delta\le Q_0$, we have
$$
\sqrt{y}-c_{22}\delta/\sqrt{Q_0}\le \sqrt{y-\delta}<
\sqrt{y+\delta} \le \sqrt{y}+c_{22}\delta/\sqrt{Q_0}.
$$
Hence,
\begin{equation}
\Pi(\delta,y)<\sum\limits_{\sqrt{y}-c_{22}\delta/\sqrt{Q_0}\le 
q \le \sqrt{y}+c_{22}\delta/\sqrt{Q_0}}
\sum\limits_{\scriptsize 
\begin{array}{cccc} (y-4\delta)rz\le m \le (y+4\delta)rz
\smallskip\\ m \equiv -bq^2\mbox{ mod } r\end{array}} 1. \label{EIN2}
\end{equation}
By Lemma 8, the double sum on the right-hand side is bounded by
\begin{equation}
\le c_{23}\sum\limits_{q\in \mathbbm{Z}} \ 
\phi\left(\frac{q-\sqrt{y}}{2c_{22}\delta/\sqrt{Q_0}} 
\right) \sum\limits_{\scriptsize \begin{array}{cccc} 
m\in \mathbbm{Z}\smallskip\\ m\equiv -bq^2 \mbox{ mod } r\end{array}} 
\phi\left(\frac{m-yrz}{8\delta rz}\right) {\rm d}y. \label{P13}
\end{equation}  
Using Lemma 7 after a linear change of variables, we transform the inner sum into 
$$
\sum\limits_{\scriptsize \begin{array}{cccc} 
m\in Z\smallskip\\ m\equiv -bq^2 \mbox{ mod } r\end{array}} 
\phi\left(\frac{m-yrz}{8\delta rz}\right) = 8\delta z
\sum\limits_{j\in \mathbbm{Z}} e\left(\frac{jbq^2}{r}+jyz\right)
\hat{\phi}(8j\delta z).
$$ 
Therefore, the double sum in (\ref{P13}) is
\begin{equation}
= 8\delta z \sum\limits_{j\in \mathbbm{Z}} e(jyz)\hat{\phi}(8j\delta z) 
\sum\limits_{d=1}^{r^*} e\left(\frac{j^*bd^2}{r^*}\right)
\sum\limits_{\scriptsize \begin{array}{cccc} 
k\in \mathbbm{Z}\\ k\equiv d \mbox{ mod } r^*\end{array}}
\phi\left(\frac{k-\sqrt{y}}{2c_{22}\delta/\sqrt{Q_0}}\right),
\label{F1} 
\end{equation}
where $r^*:=r/(r,j)$ and $j^*:=j/(r,j)$. Again using Lemma 7 after a linear 
change of variables, we transform the inner sum in (\ref{F1}) into
\begin{equation}
\sum\limits_{\scriptsize \begin{array}{cccc} 
k\in \mathbbm{Z}\\ k\equiv d \mbox{ mod } r^*\end{array}}
\phi\left(\frac{k-\sqrt{y}}{2c_{22}\delta/\sqrt{Q_0}}\right)=
\frac{2c_{22}\delta}{r^*\sqrt{Q_0}} \sum\limits_{l\in \mathbbm{Z}} 
e\left(l\cdot\frac{d-\sqrt{y}}{r^*}\right)\hat{\phi}
\left(\frac{2c_{22}l\delta}{r^*\sqrt{Q_0}}\right).\label{F2}
\end{equation}
From (\ref{F1}) and (\ref{F2}), we obtain 
\begin{eqnarray}
& & \frac{1}{\delta} \int\limits_{Q_0}^{2Q_0} \
\sum\limits_{q\in \mathbbm{Z}} \ 
\phi\left(\frac{q-\sqrt{y}}{2c_{22}\delta/\sqrt{Q_0}} 
\right) \sum\limits_{\scriptsize \begin{array}{cccc} 
m\in \mathbbm{Z}\\ m\equiv -bq^2 \mbox{ mod } r\end{array}} 
\phi\left(\frac{m-yrz}{8\delta rz}\right) {\rm d}y  \label{F3}\\
&\le& \frac{16c_{22}\delta z}{\sqrt{Q_0}} \sum\limits_{j\in \mathbbm{Z}} 
\frac{\hat{\phi}(8j\delta z)}{r^*} \sum\limits_{l\in \mathbbm{Z}}
\hat{\phi}\left(\frac{2c_{22}l\delta}{r^*\sqrt{Q_0}}\right)\left\vert 
\sum\limits_{d=1}^{r^*} 
e\left(\frac{j^*bd^2+ld}{r^*}\right) \right\vert \vert E(j,l)\vert,\nonumber
\end{eqnarray} 
where 
$$
E(j,l):=
\int\limits_{Q_0}^{2Q_0} e\left(jyz-l\cdot\frac{\sqrt{y}}{r^*}\right)\
{\rm d}y.
$$
Applying Lemmas 8 and 11, we deduce that the right-hand side of (\ref{F3}) is bounded by 
\begin{equation}
\le \frac{c_{24}\delta z}{\sqrt{Q_0}} \sum\limits_{\vert j \vert 
\le 1/(8\delta z)} \frac{1}{\sqrt{r^*}} 
\sum\limits_{\vert l\vert \le r^*\sqrt{Q_0}/(2c_{22}\delta)} 
\left\vert E(j,l)\right\vert. \label{F4}
\end{equation}

We have 
$$
E(0,0)=Q_0.
$$
If $j\not= 0$, then 
$$
\vert E(j,0) \vert \le \frac{1}{\vert j\vert z}.
$$
If $l\not=0$, then 
$$
\vert E(0,l) \vert \le \frac{c_{25}Q_0^{1/2}}{\vert l\vert}
$$
by Lemma 9 (take into account that $r^*=1$ if $j=0$).
If $j\not=0$ and $l\not=0$, then Lemma 10 yields
$$
\vert E(j,l)\vert \le 
\frac{c_{26}\sqrt{r^*}Q_0^{3/4}}{\sqrt{\vert l\vert}}.
$$
Therefore, the expression in (\ref{F4}) is bounded by
\begin{eqnarray}
\label{F5}
&\le& c_{27}\delta\left(z\sqrt{Q_0}+\frac{1}{\sqrt{Q_0}}
\sum\limits_{1\le j \le 1/(8\delta z)} \frac{1}{j\sqrt{r^*}}\right.
+z\sum\limits_{1\le l\le \sqrt{Q_0}/(2c_{22}\delta)} 
\frac{1}{l}+
\\ & & \left. zQ_0^{1/4}\sum\limits_{1\le j \le 1/(8\delta z)}\
\sum\limits_{1\le l\le r^*\sqrt{Q_0}/(2c_{22}\delta)}\frac{1}{\sqrt{l}}\right)
\nonumber\\
&\le& c_{28}\left(\delta z\sqrt{Q_0}+
\frac{\delta}{\sqrt{Q_0}}
\sum\limits_{1\le j \le 1/(8\delta z)} \frac{1}{j\sqrt{r^*}}
+\delta z\Delta^{-\varepsilon}+\right.\nonumber\\
& & \left. z\sqrt{\delta} Q_0^{1/2}
\sum\limits_{1\le j \le 1/(8\delta z)} \sqrt{r^*}\right).\nonumber
\end{eqnarray} 

Now, we evaluate the sums over $j$ in the last line of (\ref{F5}).
By the definition of $r^*$, we have 
\begin{eqnarray}
\sum\limits_{1\le j \le 1/(8\delta z)} \frac{1}{j\sqrt{r^*}}&=&
\frac{1}{\sqrt{r}}
\sum\limits_{t\vert r} \sqrt{t} \sum\limits_{\scriptsize \begin{array}{cccc} 
1\le j \le 1/(8\delta z)\smallskip
\\ (r,j)=t\end{array}} \frac{1}{j}\label{F6} \\
&\le& \frac{c_{29}\log(2+1/(8\delta z))}{\sqrt{r}}
\sum\limits_{t\vert r} \frac{1}{\sqrt{t}}\nonumber\\
&\le& c_{30}\Delta^{-\varepsilon}r^{-1/2}\nonumber
\end{eqnarray}   
and 
\begin{eqnarray}
\sum\limits_{1\le j \le 1/(8\delta z)} \sqrt{r^*}&=&
{\sqrt{r}}
\sum\limits_{t\vert r} \frac{1}{\sqrt{t}} 
\sum\limits_{\scriptsize \begin{array}{cccc} 
1\le j \le 1/(8\delta z)\smallskip\\ (r,j)=t\end{array}} 1 \label{F7} \\
&\le& \frac{\sqrt{r}}{8\delta z}
\sum\limits_{t\vert r} \frac{1}{t^{3/2}}\nonumber\\
&\le& \frac{c_{31}\sqrt{r}}{\delta z}.\nonumber
\end{eqnarray}   

Combining (\ref{EIN1}), (\ref{EIN2}), (\ref{P13}), (\ref{F3}), 
(\ref{F4}), (\ref{F5}), (\ref{F6}) and (\ref{F7}), we
obtain
\begin{equation}
P\left(\frac{b}{r}+z\right)\le 
c_{32}\Delta^{-\varepsilon}\left(1+\delta z\sqrt{Q_0}
+\delta Q_0^{-1/2}r^{-1/2}+
\delta^{-1/2}Q_0^{1/2}\sqrt{r}\right).\label{Q1}
\end{equation} 
Choosing $\delta:=Q_0\Delta/z$, we infer the desired estimate from 
(\ref{Q1}) and (\ref{P2}).$\Box$\medskip\\
 
{\bf Acknowledgement.} This paper was written during
postdoctoral stays at the Harish-Chandra Research Institute at Allahabad 
(India) and the Department of Mathematics and Statistics at Queen's University 
in Kingston (Canada). 
The author wishes to thank these institutions for financial support.

\end{document}